\documentclass[german,english,reqno]{amsart}
\usepackage{palatino}
\usepackage[T1]{fontenc}
\usepackage[latin1]{inputenc}
\usepackage{amssymb}
\IfFileExists{url.sty}{\usepackage{url}}
                      {\newcommand{\url}{\texttt}}

\makeatletter

 \theoremstyle{plain}
 \theoremstyle{plain}    
 \newtheorem{thm}{Theorem} 
 \theoremstyle{plain}    
 \newtheorem{lem}{Lemma} 
 \theoremstyle{remark}
 \newtheorem*{rem*}{Remark}
 \theoremstyle{plain}    
 \newtheorem*{cor*}{Corollary}

\usepackage{color}
\usepackage[marginal]{footmisc}
\usepackage{varioref}
\usepackage{array}

\newcommand{\ZZ}{\mathbb Z}

\DeclareMathOperator{\real}{Re}

\DeclareMathOperator{\spann}{span}
\DeclareMathOperator{\vol}{vol}
\DeclareMathOperator{\GL}{GL}


\usepackage{babel}
\makeatother
\begin{document}

\title{An asymptotic expansion for the discrete harmonic potential}

\author{Gady Kozma}

\thanks{This work is part of the research program of the European Network
{}``Analysis and Operators'', contract HPRN-CT-00116-2000 supported
by the European Commission.}

\email{gady@post.tau.ac.il, gadykozma@hotmail.com}

\curraddr{Tel Aviv University, Tel Aviv, Israel.}

\author{Ehud Schreiber}

\email{schreib@physics.ubc.ca}\address{Tel Aviv University, Tel Aviv, Israel.}

\begin{abstract}
We give two algorithms that allow to get arbitrary precision asymptotics
for the harmonic potential of a random walk.
\end{abstract}


\maketitle
\theoremstyle{plain}
\newtheorem*{question}{Question}

\section{Introduction}

The discrete harmonic potential $a$ of a random walk $R$, starting
from $0$, on a lattice $Z$ can be most easily defined using \begin{equation}
a(z):=\sum_{n=0}^{\infty}P_{n}(0)-P_{n}(z)\label{eq:defa}\end{equation}
where $P_{n}(z)$ is the probability of $R$ to be at $z$ at the
$n$th step. It is easy to conclude from (\ref{eq:defa}) that $\Delta a=\delta_{\{0\}}$
where $\Delta$ is the discrete Laplacian for the walk $R$. Less
obvious are the facts that $a$ is positive, and that these two properties
determine $a$ up to the addition of a positive constant. {}``Positive''
may be replaced by {}``with sub-linear growth'' in dimension greater
than $1$. 

The discrete harmonic potential is an interesting quantity, strongly
related to the discrete Green function $g(x,y)$: on the one hand,
$a(x)=g(x,x)$ where $g$ is the Green function corresponding to the
set $\{0\}$ \cite[P11.6, page 118]{S76}. On the other hand, $g$
for arbitrary sets can be calculated if $a$ is known --- when the
set is finite there is an explicit formula, \cite[T14.2, page 143]{S76}.

Thus it becomes interesting to calculate or estimate $a$. Except
in dimension $1$ where $a(x)=|x|$, the problem is far from trivial.
Take as an example the harmonic potential of the regular random walk
on $\ZZ^{2}$ which has a logarithmic nature. The estimate of $a$
most common in the literature is \begin{equation}
a(z)=\frac{2}{\pi}\log|z|+\lambda+O(|z|^{-2})\label{eq:stoher}\end{equation}
where $\lambda$ is some number, which can be expressed using Euler's
$\gamma$ constant, $\lambda=\frac{2}{\pi}\gamma+\frac{1}{\pi}\log8$.
The earliest proof of (\ref{eq:stoher}) we are aware of is in St\"ohr
\cite{S49}. A more accessible proof of a weaker result, giving an
error estimate of $o(1)$, can be found in Spitzer \cite[P12.3, page 124]{S76}.
Spitzer's proof (which is not specific to the lattice $\ZZ^{2}$),
like St\"ohr's, relies on the fact that the Fourier transform $\widehat{P_{n}}$
has an explicit formula, which can be summed to give a {}``pseudo-Fourier''
representation of $a$ as the explicit integral\begin{equation}
a(z)=\int\frac{1-e^{i\langle\xi,z\rangle}}{1-\widehat{P_{1}}(\xi)}\, d\xi\quad;\label{eq:psdfur}\end{equation}
 indeed, this is how he proves that the sum in (\ref{eq:defa}) converges
in the first place. Fukai and Uchiyama used this technique to find
the next order approximation --- in the case of the simple random
walk on $\mathbb{Z}^{2}$ it is $-\frac{\real z^{4}}{6\pi|z|^{6}}$
--- and to show that an asymptotic polynomial expansion exists. See
\cite{FU96} for the two dimensional case and \cite{U98} for the
case $d\geq3$.

The purpose of this note is to give two alternative approaches to
the computation of $a$, which allow to get asymptotics of arbitrary
precision. For example, here are the first few terms for the regular
random walk on $\mathbb{Z}^{2}$:\begin{eqnarray}
a(z) & = & \frac{2}{\pi}\log|z|+\lambda-\real\left(\frac{z^{4}}{6\pi|z|^{6}}\:+\right.\nonumber \\
 &  & +\:\frac{3z^{4}}{20\pi|z|^{8}}+\frac{5z^{8}}{24\pi|z|^{12}}+\nonumber \\
 &  & +\:\frac{51z^{8}}{56\pi|z|^{14}}+\frac{35z^{12}}{36\pi|z|^{18}}+\nonumber \\
 &  & +\:\frac{217z^{8}}{160\pi|z|^{16}}+\frac{45z^{12}}{4\pi|z|^{20}}+\frac{1925z^{16}}{192\pi|z|^{24}}+\nonumber \\
 &  & +\: O(|z|^{-10})\bigg)\quad.\label{eq:a_high_order}\end{eqnarray}
The first approach is a brute-force one which starts from (\ref{eq:defa}):
$P_{n}$ may be represented as a sum of multinom coefficients, and
plugging in Stirling's expansion, and some elementary algebra gives
the representation (\ref{eq:a_high_order}). Unfortunately, the algebra
involved is too long. In effect, for all but the first coefficient
(the $-\frac{\real z^{4}}{6\pi|z|^{6}}$), it is simply not practical
to do the calculations by hand. We demonstrate this approach in two
ways: we use it to prove that a polynomial approximation exists, reproducing
the results of \cite{FU96,U98} in the case the random walk is bounded
(Fukai and Uchiyama show these results under the weaker assumption
that each step of the random walk has $K+d+\epsilon$ moments, where
$K$ is the required precision in the expansion, $d$ is the dimension
and $\epsilon>0$). See theorem \ref{thm:existpoly} on page \pageref{thm:existpoly};
and we use it to get some explicit constants and high order expansions
for other walks. See section \ref{sec:Explicit-constants}.

The second approach is to approximate the discrete Laplacian $\Delta$
using an appropriate differential operator. For example, for the regular
random walk on $\mathbb{Z}^{2}$ we have, from the Taylor expansion,\begin{eqnarray}
(\Delta f)(z) & = & {\textstyle \frac{1}{4}}\big(f(z+1)+f(z-1)+f(z+i)+f(z-i)\big)-f(z)\nonumber \\
 & = & \frac{1}{4}\left(\frac{\partial^{2}f}{\partial x^{2}}+\frac{\partial^{2}f}{\partial y^{2}}\right)+\frac{1}{48}\left(\frac{\partial^{4}f}{\partial x^{4}}+\frac{\partial^{4}f}{\partial y^{4}}\right)+\dotsb\label{eq:disc_dif}\end{eqnarray}
In fact, $\Delta=\frac{1}{2}(\cosh\frac{\partial}{\partial x}+\cosh\frac{\partial}{\partial y})-\mathbf{I}$,
but we will not use this representation. Given that the expansion
of $a$ has an appropriate polynomial form, the expansion (\ref{eq:disc_dif})
allows to get many equations on the coefficients in (\ref{eq:a_high_order}).
These equations are of course insufficient --- the equation $\Delta a=\delta_{\{0\}}$
does not determine $a$ uniquely, and the condition that $a$ is positive
or slowly increasing is hard to encode into the differential equations.
However these differential equations allow to prove many of the apparent
properties of (\ref{eq:a_high_order}), most notably that the coefficient
of $z^{l}|z|^{-k}$ for $l<\frac{1}{2}k$ is always zero. This is
theorem \ref{thm:klhalf}. Unfortunately, it is still not clear why
all the coefficients are negative.

Returning to the case of the regular random walk on $\mathbb{Z}^{2}$,
it can be seen that the equation $\Delta a=\delta_{\{0\}}$ combined
with the $\frac{\pi}{2}$-rotational symmetry gives that the values
of $a$ are uniquely defined by their values on one diagonal, say
$\{ m+im\}_{m=0}^{\infty}$, and vice versa, any arbitrary sequence
on the diagonal can be extended to a symmetric function $f$ with
$\Delta f=\delta_{\{0\}}$ using a very simple recursion. On the other
hand, it turns out that the values of $a$ on the diagonal can be
calculated explicitly using the integral (\ref{eq:psdfur}) to get

\begin{equation}
a(m+im)=\frac{4}{\pi}\left(1+\frac{1}{3}+\frac{1}{5}+\dotsb+\frac{1}{2m-1}\right)\label{eq:anexplic}\end{equation}
 \cite[chapter 15]{S76}. Indeed we have just sketched an algorithm
for calculating any value of $a(z)$: this is the McCrea-Whipple algorithm
\cite{MW40}%
\footnote{Can be found in Spitzer as well: \cite[chapter 15, page 148]{S76}%
}. As is now evident, all the values are of the form $n+\frac{1}{\pi}q$,
$n\in\mathbb{Z}$ and $q\in\mathbb{Q}$. Interestingly, when not on
the diagonal, $n$ and $q$ increase exponentially, even though the
result is only logarithmic in size%
\footnote{No, this is not a very efficient way to calculate $\pi$, it requires
at least $k^{3}$ operations to get a precision of $k$ digits...%
}. As for the coefficients of (\ref{eq:a_high_order}), it turns out
--- perhaps unsurprisingly --- that (\ref{eq:anexplic}) can be used
to complete the few coefficients that cannot be deduced from (\ref{eq:disc_dif}).
This is theorem \ref{thm:eqssuff}.

\label{Thirdapproach}There is a third approach to the problem which
we shall not discuss at length as it seems inferior to both former
ones. Roughly it goes as follows (for the case of the regular random
walk on $\mathbb{Z}^{2}$). {}``Guess'' that the solution is approximately
$\frac{2}{\pi}\log|z|$. We discretize to $\mathbb{Z}^{2}$ in the
natural way: define \[
f(z):=\begin{cases}
\frac{2}{\pi}\log|z| & z\neq0\\
-1 & z=0\end{cases}\quad.\]
 Actually, it is not necessary to guess the exact value $\frac{2}{\pi}$,
other constants not too different also work. A calculation can show
that \begin{equation}
\sum_{z\in\mathbb{Z}^{2}}|(\Delta f-\delta_{\{0\}})(z)|<1\quad.\label{eq:betalt1}\end{equation}
This allows us to write a series of corrections of $f$ as follows:
$f_{1}=f$ and then \[
f_{n}=f_{n-1}-(\Delta f_{n-1}-\delta_{\{0\}})*f\]
where $*$ is the convolution on $\mathbb{Z}^{2}$. Since $\Delta(g*h)=g*\Delta h$
we get \[
\Delta f_{n}-\delta_{\{0\}}=(\Delta f_{n-1}-\delta_{\{0\}})*(\Delta f-\delta_{\{0\}})\]
and hence (\ref{eq:betalt1}) gives that the $f_{n}$'s converge%
\footnote{The convergence rate is obviously exponential, but one might also
define $f_{n}=f_{n-1}-(\Delta f_{n-1}-\delta_{\{0\}})*f_{n-1}$ and
get double exponential convergence.%
} in the $L^{1}$ norm to an $a$ which satisfies $\Delta a=\delta_{\{0\}}$
and $|a|\leq C\log|z|+C$ so it is the harmonic potential up to an
additive constant. Furthermore, the fact that $\Delta f-\delta_{\{0\}}=O(|z|^{-4})$
can be used to derive similar estimates for $a$. However, the best
approximation we were able to get from that method is \begin{equation}
a(z)=\frac{2}{\pi}\log|z|+\lambda-\real\frac{z^{4}}{6\pi|z|^{6}}+O(|z|^{-4}\log|z|)\label{eq:thirdstoher}\end{equation}
and getting this approximation was not significantly easier than the
brute force approach. Moreover, this approach did not allow to get
actual values for any of the constants --- only that some constants
exist, and their values were derived from (\ref{eq:anexplic}).

We have described the contents of this paper except for the last section.
That section contains the results of some computer-aided simulations
and a calculation of an \emph{explicit} constant in the $O(\cdot)$
in (\ref{eq:stoher}).

We wish to thank Greg Lawler for some useful remarks.

\subsection{Standard definitions}

A lattice is a discrete additive subgroup of $\mathbb{R}^{n}$. The
dimension of a lattice $Z$ (denoted by $\dim Z$) is the dimension
of the linear span of $Z$ as a linear subspace of $\mathbb{R}^{n}$.
A basis for a $d$-dimensional lattice is a set $\{ e_{1},\dotsc,e_{d}\}$
such that $Z=e_{1}\mathbb{Z}+\dotsb+e_{d}\mathbb{Z}$. $Z'$ is called
a sublattice of $Z$ if it is a subgroup of $Z$, and its index is
the size of $Z/Z'$. The volume of (the cell of) a lattice $Z$, denoted
by $\vol Z$, is the volume of $\mathbb{R}^{d}/Z$ (the discreteness
of $Z$ allows us to select a measurable set of representatives).
Alternatively it can be defined as \[
\lim_{r\rightarrow\infty}\frac{|B(r,0)|}{\#(B(r,0)\cap Z)}\]
where $B(r,0)$ is a ball of radius $r$ around $0$ and $|B(r,0)|$
is its volume.

A random walk $R$ on a lattice $Z$ is a probabilistic process $\{ R_{0},R_{1},\dotsc\}$,
$R_{i}\in Z$, such that $R_{0}=0$ and $R_{i}-R_{i-1}$ has a distribution
independent of $i$, and such that for every $z\in Z$ there exists
some $n$ such that $\mathbb{P}(z=R_{n})>0$.%
\footnote{This definition is a bit restrictive. For example, a $1$-dimensional
random walk which always goes to the right is not a {}``random walk
on a lattice''.%
} It is bounded if $R_{i}-R_{i-1}$ is bounded. We will typically confuse
the process $R$ with the distribution of any step (e.g.~with the
distribution of $R_{1}$) so that we can use notations such as $\mathbb{E}R$
comfortably. The dimension of $R$ is the dimension of the lattice
$Z$, and is denoted by $\dim R$.

The drift of a random walk $R$ is $\mathbb{E}R$. The random walk
is balanced if its drift is zero. 

The discrete Laplacian $\Delta$ of a random walk $R$ is an operator
on functions on the lattice $Z$ defined by \[
(\Delta f)(z)=\mathbb{E}f(z+R_{1})-f(z)\quad.\]
Functions with $\Delta f\equiv0$ are called (discretely) harmonic.

The continuous Laplacian on $\mathbb{R}^{d}$, $\frac{\partial^{2}}{\partial z_{1}^{2}}+\dotsb+\frac{\partial^{2}}{\partial z_{d}^{2}}$
will be denoted by $\Delta_{C}$. Functions with $\Delta_{C}f\equiv0$
will be called continuously harmonic.

$||v||$ will always refer to the $L^{2}$ norm. Similarly, for a
matrix $A$, $||A||$ refers to its norm as an operator from $L^{2}$
to $L^{2}$.

\section{The direct approach}

The aim of this section is to prove the following theorem.

\begin{thm}
\label{thm:existpoly}Let $a$ be the harmonic potential of a $d$-dimensional
balanced bounded random walk on a lattice $Z$. Let $K>0$ be some
integer. Then there exists a constant $\tau$, a linear map $A$,
a polynomial $Q=Q_{K}$ in $d$ variables and an integer $L=L_{K}$
such that\begin{equation}
a(z)=\tau\mu_{d}(||Az||)+\frac{Q(z)}{||Az||^{L}}+O(||z||^{-K})\label{eq:thm1}\end{equation}
where \[
\mu_{d}(z):=\begin{cases}
||z||^{2-d} & d\neq2\\
\log||z|| & d=2\end{cases}\quad.\]

\end{thm}
$\tau$ and $A$ have explicit formulas --- see the comments after
the end of the proof (page \pageref{constants} below). The expression
$Q(z)/||Az||^{L}$ is $\lambda+O(||Az||^{1-d})$ (this is classical,
of course, but will also follow from the proof below). The number
$\lambda$ is of course also walk-specific.

In the case $d=1$ the equation $\Delta a=\delta_{\{0\}}$ deteriorates
into a simple recursion formula and the theorem is nothing but an
exercise (and $Q\equiv\lambda$). Thus we concentrate on the case
$d>1$.

\begin{lem}
\label{lem:binom}Let $p_{1},\dotsc,p_{k}\in\left]0,1\right[$ such
that $\sum p_{i}=1$. Let $K\in\mathbb{N}$. Then there exists a $Q=Q_{p,K}$
such that\begin{align}
\lefteqn{B_{p}(n,p_{1}n+w_{1},\dotsc,p_{k}n+w_{k})=}\label{eq:binom}\\
 & \qquad\qquad\frac{1}{(2\pi n)^{(k-1)/2}\sqrt{p_{1}\dotsb p_{k}}}e^{-\sum w_{i}^{2}/2p_{i}n}Q(n,w_{1},\dotsc,w_{k})+O(n^{-K})\nonumber \end{align}
for every $w_{1},\dotsc,w_{k}$ such that $\sum w_{i}=0$. $B$ here
is the standard multinom\[
B_{p}(n,m_{1},\dotsc,m_{k})=\frac{n!}{m_{1}!\cdot\dotsb m_{k}!}p_{1}^{m_{1}}\cdot\dotsb\cdot p_{k}^{m_{k}},\quad m_{1}+\dotsb+m_{k}=n\]
and $Q$ satisfies that $n^{L}Q$ is a polynomial for some $L=L_{K}$
and $Q(n,w)=1+o(1)$ uniformly in $||w||\leq C\sqrt{n\log n}$.
\end{lem}
The proof of theorem \ref{thm:existpoly} is replete with these rational
functions $Q$ obeying the condition above. Notice that this simply
means that \[
Q(n,w_{1},\dotsc,w_{k})=\sum_{i,\vec{j}}c_{i\vec{j}}n^{i}w_{1}^{j_{1}}\dotsm w_{k}^{j_{k}}\]
such that $c_{0\vec{0}}=1$ and otherwise $c_{i\vec{j}}\neq0$ only
if $-L\leq i<0$ and $0\leq j_{1}+\dotsb+j_{k}<-2i$. It turns out
(and this can be deduced from the proofs below with some care) that
$L_{K}=4K$ is always sufficient. However we will have no use for
this fact.

\begin{proof}
[Proof of lemma \ref{lem:binom}]Consider Stirling's series,\[
n!=\frac{n^{n}}{e^{n}}\sqrt{2\pi n}\cdot\exp\left(Q_{1}(n)+O(n^{-K})\right)\]
where $n^{K-1}Q_{1}$ is some polynomial with $Q_{1}(n)=o(1)$ as
$n\rightarrow\infty$. This gives\begin{align*}
\lefteqn{B_{p}(n,p_{1}n_{1}+w_{1},\dotsc,p_{k}n_{k}+w_{k})=\frac{1}{(2\pi n)^{(k-1)/2}\sqrt{p_{1}\dotsb p_{k}}}\prod_{i=1}^{k}\frac{1}{\sqrt{1+\frac{w_{i}}{p_{i}n}}}\;\cdot}\\
 & \qquad\exp\left(-\sum_{i=1}^{k}(p_{i}n+w_{i})\log\left(1+\frac{w_{i}}{p_{i}n}\right)+Q_{1}(n)-\sum_{i=1}^{k}Q_{1}(p_{i}n+w_{i})+O(n^{-K})\right)\end{align*}
 Assume for the moment that $||w||<C\sqrt{n\log n}$ where $C$ will
be fixed later. In this case we can insert monomials such as $\frac{w_{i}^{2K+1}}{n^{2K+1}}$
into $O(n^{-K})$. Thus, each of the factors $\frac{1}{\sqrt{1+\frac{w_{i}}{p_{i}n}}}$
equals $Q_{2}(n,w_{i})+O(n^{-K})$ with $n^{2K}Q_{2}$ some polynomial
with $Q_{2}=1+o(1)$ ($Q_{2}$ depends on $p_{i}$, of course). The
term inside the exponent evaluates to $-\sum\frac{w_{i}^{2}}{2p_{i}n}+Q_{3}(n,\vec{w})+O(n^{-K})$
with $n^{2K}Q_{3}$ a polynomial satisfying $Q_{3}=o(1)$ uniformly
in $||w||<C\sqrt{n\log n}$ as $n\rightarrow\infty$. This proves
the case $||w||<C\sqrt{n\log n}$. In the case $||w||>C\sqrt{n\log n}$
both sides of (\ref{eq:binom}) are $O(n^{-K})$ which follows easily
from the fact that in this case, for $C$ sufficiently large\[
\exp\left(-\sum\frac{w_{i}^{2}}{2p_{i}n}\right)\leq\exp\left(-C\log n\right)\leq Cn^{-K}\quad.\qedhere\]

\end{proof}
\begin{lem}
\label{lem:sumexp}Let $v\in\mathbb{R}^{d}$ and $A\in\GL(\mathbb{R}^{d})$.
Then\[
\sum_{z\in\mathbb{Z}^{d}}\exp\!\left(-\frac{1}{n}||A(z+v)||^{2}\right)=\frac{(\pi n)^{d/2}}{|A|}+O(e^{-cn})\]
where the $O(\cdot)$ is uniform in $v$. ($|A|$ here is the determinant
of $A$). 

Furthermore, if $P$ is any polynomial and $y\in\mathbb{R}^{d}$ then\begin{equation}
\sum_{z\in\mathbb{Z}^{d}+y}\exp\!\left(-\frac{1}{n}||A(z+v)||^{2}\right)P(z)=n^{d/2}Q(n)+O(e^{-cn}||v||^{\deg P})\label{eq:withP}\end{equation}
where $Q$ is a polynomial. The coefficients of $Q$ are fixed polynomial
functions (depending on $A$) of the coefficients of $P$ and of $v$.
$Q$ is independent of $y$ and the $O(\cdot)$ is uniform in $v$
and $y$.
\end{lem}
The use of the parameter $n$ in the formulation of the lemma is a
bit artificial. Later, when we shall use the lemma, $n$ will be an
integer number (the number of steps), but this fact is not needed
for the formulation or proof of the lemma.

\begin{proof}
We use the $d$-dimensional Poisson summation formula, $\sum f(\vec{n})=\sum\widehat{f}(2\pi\vec{n})$
(the Fourier transform being defined by $\widehat{f}(\xi)=\int f(x)e^{i\langle x,\xi\rangle}\, dx$)
for $f(z)=\exp\left(-\frac{1}{n}||A(z+v)||^{2}\right)(z+v)^{\vec{r}}$
where $\vec{r}=(r_{1},\dotsc,r_{d})\in\{0,1,\dotsc\}^{d}$ and $(z+v)^{\vec{r}}$
is short for $(z_{1}+v_{1})^{r_{1}}\dotsm(z_{d}+v_{d})^{r_{d}}$.
A simple calculation shows that\begin{equation}
\widehat{f}(\xi)=e^{-i\langle v,\xi\rangle}i^{-r}\frac{d^{r}}{d\xi^{\vec{r}}}\frac{(\pi n)^{d/2}}{|A|}\exp\!\left(-{\textstyle \frac{1}{4}}n\left\Vert \left(A^{-1}\right)^{*}\xi\right\Vert ^{2}\right)\label{eq:fhat}\end{equation}
where $r=r_{1}+\dotsb+r_{d}$, and as usual $\frac{d^{r}}{d\xi^{\vec{r}}}=\frac{d^{r_{1}}}{d\xi_{1}^{r_{1}}}\dotsb\frac{d^{r_{d}}}{d\xi_{d}^{r_{d}}}$.
This means that as $n\rightarrow\infty$ we have \[
\sum_{\xi\in2\pi\mathbb{Z}^{d}}\widehat{f}(\xi)=\widehat{f}(0)+O\!\left(n^{r+d/2}\left\Vert A^{-1}\right\Vert ^{d+2r}e^{-n/4||A||}\right)\]
(note that $|A|^{-1}\leq||A^{-1}||^{d}$). The exponent inside the
$O(\cdot)$ kills of course all the other factors. The case $\vec{r}=\vec{0}$
immediately gives the first part of the lemma. For the second part,
put $\xi=0$ in (\ref{eq:fhat}) and get that $\widehat{f}(0)\neq0$
only when $r$ is even, and that $\widehat{f}(0)=n^{d/2}Q_{\vec{r}}(n)$
where $Q_{\vec{r}}$ depends only on $A$. Writing $P=\sum_{\vec{r}}c_{\vec{r}}(z+v)^{\vec{r}}$
we get $Q=\sum_{\vec{r}}c_{\vec{r}}Q_{\vec{r}}$ which concludes the
lemma: the fact that $c_{\vec{r}}$ are polynomial in the coefficients
of $P$ and in $v$ gives the same for the coefficients of $Q$; the
fact that $c_{\vec{r}}=O(||v||^{\deg P})$ gives the same in the $O(\cdot)$
in (\ref{eq:withP}). The result does not depend on $y$ (except in
the $O(\cdot)$ error) since $\sum f(n)$ does not depend on $v$
--- in effect the dependence on $v$ appears only via the dependence
of $c_{\vec{r}}$ on $v$.
\end{proof}
\begin{lem}
\label{lem:multidim}Let $R$ be a $d$-dimensional balanced bounded
random walk on the lattice $Z$. Let $K\in\mathbb{N}$. Then there
exists a sublattice $Y\subset\mathbb{Z}\times Z$, a constant $\tau'$,
a quadratic function $B$ and a $Q$ such that \begin{equation}
P_{n}(z)=\begin{cases}
\tau'n^{-d/2}e^{-B(z)/n}Q(n,z)+O(n^{-K}) & (n,z)\in Y\\
0 & (n,z)\not\in Y\end{cases}\quad.\label{eq:reslem}\end{equation}
with $n^{L}Q$ a polynomial in $n$ and in the coordinates of $z$,
and $Q(n,z)=1+o(1)$ uniformly in $||z||\leq C\sqrt{n\log n}$.
\end{lem}
Assume for simplicity that $Z\subset\mathbb{R}^{d}$. A simple comparison
of the expectations of the left and right side of (\ref{eq:reslem})
shows that $B(z)=\frac{1}{2}\langle M^{-1}z,z\rangle$, where $M$
is the correlation matrix of $R$, $M_{ij}=\mathbb{E}\langle R,e_{i}\rangle\langle R,e_{j}\rangle$,
where $e_{i}$ are unit vectors. Notice that the expression $\langle M^{-1}z,z\rangle$
is independent of the coordinate system chosen.

The lattice $Y$ is clearly $d+1$ dimensional since otherwise we
would have $\dim R<d$. Thus $Y\cap Z$ is $d$-dimensional and $\vol Y\cap Z$
is finite. This, with the requirement $\sum P_{n}=1$ allows to calculate
$\tau'$ and get \begin{equation}
\tau'=\frac{\vol(Y\cap Z)}{(2\pi)^{d/2}|\det M|^{1/2}}\quad.\label{eq:Cprob}\end{equation}
We will not be using these equalities in the proof, though.

\begin{proof}
By applying a linear transformation we may assume $R$ is a random
walk on $\mathbb{Z}^{d}$. Assume $R$ moves from $z$ to $z+v_{i}$
with probability $p_{i}$ for $i=1,\dotsc,N+1$. Let \[
W=\left\{ (n,w_{1},\dotsc,w_{N+1}):p_{i}n+w_{i}\in\mathbb{Z}\;\forall i\textrm{ and }\sum w_{i}=0\right\} \]
which is an $(N+1)$-dimensional lattice. For an $n\in\mathbb{N}$
let $k_{i}$ be the number of times (between $0$ and $n-1$) that
$R$ moved from $z$ to $z+v_{i}$. We use lemma \ref{lem:binom}
with some $K_{2}$ which will be fixed later. This gives us\begin{equation}
\mathbb{P}(k_{i}=np_{i}+w_{i})=Cn^{-N/2}\exp\left(-\frac{1}{n}B_{2}(w)\right)Q_{2}(n,w_{1},\dotsc,w_{N}\big)+O(n^{-K_{2}})\label{eq:Pmulti}\end{equation}
where $B_{2}$ is a strictly positive quadratic form, $C$ is some
constant and $n^{L}Q_{2}$ is a polynomial with $Q_{2}=1+o(1)$ uniformly
for $||w||\leq C\sqrt{n\log n}$. Now, for a point $z\in\mathbb{R}^{d}$
we can write\begin{align}
\mathbb{P}(R_{n}=z) & =\sum_{w\in W(n,z)}\mathbb{P}(k_{i}=np_{i}+w_{i})\label{eq:sumsublattice}\\
W(n,z) & :=\left\{ w:(n,w)\in W\;\wedge\;\sum w_{i}v_{i}=z\right\} \quad.\nonumber \end{align}
Since $W(n,z)$ are intersections of a lattice with a translated subspace,
we get that $W(n,z)=\emptyset$ for $(n,z)$ outside some lattice
$Y$, and for $(n,z)\in Y$ we get that $W(n,z)$ are translations
of one fixed lattice $X$. Since we know $Y$ is $d+1$ dimensional
(see the comment just after the statement of this lemma) we get that
$X$ is $(N-d)$-dimensional. The calculation of the sum in (\ref{eq:sumsublattice})
will be done of course with lemma \ref{lem:sumexp} but we need some
preparations first. By applying a linear transformation to the $N$-dimensional
space of $w$'s we may assume that \[
W(n,z)=z_{1}e_{1}+\dotsb+z_{d}e_{d}+e_{d+1}\mathbb{Z}+\dotsb+e_{N}\mathbb{Z}+y_{n}\]
where $y_{n}\in\spann\{ e_{d+1},\dotsc,e_{N}\}$. Combining (\ref{eq:Pmulti}),
(\ref{eq:sumsublattice}) and the linear transformation we get \begin{equation}
\mathbb{P}(R_{n}=z)=\sum_{w\in W(n,z)}Cn^{-N/2}\exp\!\left(-\frac{1}{n}B_{3}(w)\right)Q_{3}(n,w)+O(n^{-K_{3}})\label{eq:afterlin}\end{equation}
where $K_{3}=K_{2}-N+d$. This last error estimate follows easily
from the following obvious estimates:\[
\sum_{\substack{w\in W(n,z)\\
p_{i}n+w_{i}\in[0,n]}
}O(n^{-K_{2}})=O(n^{-K_{3}}),\qquad\sum_{\substack{w\in W(n,z)\\
p_{i}n+w_{i}\not\in[0,n]}
}e^{-B_{3}(w)/n}Q_{3}(n,w)=O(e^{-cn})\]

From now on it would be easier to replace $B_{3}(w)$ with $||Aw||^{2}$
for some $A\in\GL(\mathbb{R}^{N})$ which is possible since $B_{3}$
is also strictly positive. We denote $x=(0,\dotsc,0,\linebreak[0]w_{d+1},\dotsc,w_{N})$
so that $w=x+z$ and we can write\[
||Aw||^{2}=||Ax||^{2}+||Az||^{2}+2\langle x,A^{*}Az\rangle\quad.\]
Since $A^{*}A$ is strictly positive we get that its lower right $N-d$
minor is also strictly positive and thus invertible. Therefore there
exists some $v=(0,\dotsc,0,v_{d+1},\dotsc,v_{N})$ such that \[
(A^{*}Av)_{i}=(A^{*}Az)_{i},\quad i=d+1,\dotsc,N\quad.\]
This of course implies $\langle x,A^{*}Av\rangle=\langle x,A^{*}Az\rangle$
and then \begin{equation}
||Aw||^{2}=||Az||^{2}-||Av||^{2}+||A(x+v)||^{2}\quad.\label{eq:Aparts}\end{equation}
The importance of (\ref{eq:Aparts}) is that its left part, $||Az||^{2}-||Av||^{2}$,
is constant for $w\in W(n,z)$, while $\exp\left(-\frac{1}{n}||A(x+v)||^{2}\right)$
can be summed by using lemma \ref{lem:sumexp} on the last $N-d$
coordinates. Notice that it depends on the $z$'s only via $v$, which
is linear in $z$. In other words, we can now define our quadratic
form $B$:\[
B(z):=||Az||^{2}-||Av(z)||^{2}\]
 and apply lemma \ref{lem:sumexp}. We get \begin{align*}
(\ref{eq:afterlin}) & =Cn^{-N/2}e^{-B(z)/n}\sum_{x\in\mathbb{Z}^{N-d}+y_{n}}\exp\!\left(-\frac{1}{n}||A(x+v)||^{2}\right)Q_{3}(n,z,x)+O(n^{-K_{3}})\\
 & =Cn^{-d/2}e^{-B(z)/n}Q(n,z)+O(n^{-K_{3}})+O(e^{-cn}||v||^{\deg Q_{3}})\end{align*}
We may already remark that picking $K_{2}=K+N-d$ will give us a precision
of $O(n^{-K})$ as required. Obviously $O(e^{-cn}||v||^{\deg Q_{3}})=O(n^{-K})$
for any $K$ as $v$ is linear in $z$ and $||z||\leq C\sqrt{n\log n}$.
To say something about $Q$, we need to consider $Q_{3}$ as a polynomial
in $x$ with coefficients which are polynomials in $z$ and negative
powers in $n$, i.e.\[
Q_{3}=\sum_{\vec{r}}c_{\vec{r}}(n,z)x^{\vec{r}}\]
where $n^{L}c_{\vec{r}}$ are polynomials. Lemma \ref{lem:sumexp}
tells us that $n^{L}Q$ is a polynomial in $n$ and that its coefficients
are polynomial in the coefficients of $Q_{3}$ and in $v$ (which
is linear in $z$). Thus we get that $n^{L}Q$ is polynomial in both
$n$ and $z$. This concludes the lemma: the only claim we haven't
shown is that $Q=1+o(1)$ which can be deduced, for example, from
the central limit theorem%
\footnote{Actually, it is also possible to deduce it from the arguments above
with a little care as to the exact dependence between $P$, $Q$ and
$v$ in lemma \ref{lem:sumexp}.%
}.
\end{proof}
\begin{rem*}
For specific walks it is possible to prove lemma \ref{lem:multidim}
using much simpler techniques. For example, for the simple random
walk on $\mathbb{Z}^{2}$, we have \[
P_{n}(x+iy)=4^{-n}\sum_{k=0}^{n}{{n \choose k}}{{k \choose (k-x)/2}}{{n-k \choose (n-k-y)/2}}\]
where we make the notational convention that ${{n \choose \alpha}}=0$
if $\alpha$ is not an integer. This formula allows to get lemma \ref{lem:multidim}
from lemma \ref{lem:binom} more-or-less directly, with no need to
go through lemma \ref{lem:sumexp}. This works for other walks as
well. However it does not seem to work for the very interesting case
of the random $6$-walk on the triangular lattice, a walk which can
be considered no less natural than the random walk on $\mathbb{Z}^{2}$.
We mean here a two dimensional random walk that goes at each step
with probability $\frac{1}{6}$ from $z$ to one of the points $z+\omega^{i}$,
$i=0,\dotsc,5$ where $\omega=\frac{1}{2}+i\frac{\sqrt{3}}{2}=\sqrt[6]{1}$.
\end{rem*}
\begin{lem}
\label{lem:skhumnae}Let $q\in\mathbb{N}$ and $1\leq j\leq q$. For
every $s>1$ and for $r\rightarrow\infty$ we have the superpolynomial
estimate\begin{equation}
F_{s,j,q}(r):=\sum_{\substack{n\equiv j\;(q)\\
n\geq j}
}n^{-s}e^{-r/n}=\frac{\Gamma(s-1)}{qr^{s-1}}+O(r^{-K})\quad\forall K;\label{eq:Fsgt1}\end{equation}
while for $s=1$ we have\[
F_{1,j,q}(r):=\sum_{\substack{n\equiv j\;(q)\\
n\geq j}
}\frac{1}{n}\left(e^{-r/n}-1\right)=q^{-1}\log r+C+O(r^{-K})\quad\forall K.\]

\end{lem}
\begin{proof}
Again we use the (one dimensional) Poisson summation formula, this
time for the function \[
f(x)=\begin{cases}
x^{-s}e^{-1/x} & x>0\\
0 & x\leq0\end{cases}\quad.\]
Simple Fourier operations show that\[
F_{s,j,q}=\frac{1}{qr^{s-1}}\sum_{\xi\in\mathbb{Z}}e^{-2\pi ij\xi/q}\widehat{f}\!\left(2\pi\frac{r\xi}{q}\right)\quad.\]
Since $f\in C^{\infty}$ we get that $\widehat{f}(\xi)=O(\xi^{-K})$
for all $K$ and hence $F=\frac{r^{1-s}}{q}\widehat{f}(0)+O(r^{-K})$.
This concludes the case $s>1$. For the case $s=1$ we notice that
the $O(\cdot)$ in (\ref{eq:Fsgt1}) is uniform in $s\in(1,2]$ since
it depends only on the $L^{1}$ norm of the $K$'th derivative of
$f$. Thus we need only to calculate \[
\lim_{s\rightarrow1}\frac{\Gamma(s-1)}{r^{s-1}}-\zeta(s)\]
where $\zeta$ is Riemann's function. Writing $\Gamma=\frac{1}{s}+C+o(1)$
and $\zeta(s)=\frac{1}{s-1}+C+o(1)$ and a little use of l'H\^{o}pital
rule will give $\log r+C$ and conclude the lemma.
\end{proof}
\begin{rem*}
The exact asymptotics of the error are harder to figure out. This
is of course a question about the asymptotics of $\widehat{f}$. A
much better estimate than $\widehat{f}(\xi)\leq\xi^{-K}$ can be had
by changing the path of integration of $x^{-s}e^{-1/x+ix\xi}$ from
$[0,\infty)$ to $[0,1+i]\cup[1+i,\infty+i)$ or to its conjugate
(depending on the sign of $\xi$). On the segment $[0,\xi^{-1/2}(1+i)]$
we have $e^{-1/x}=O(e^{-c\sqrt{\xi}})$. On the reminder we have $e^{ix\xi}=O(e^{-c\sqrt{\xi}})$.
So we get $\widehat{f}(\xi)=O(e^{-c\sqrt{\xi}})$.

In the other direction, there are known precise results connecting
the {}``depth'' of the zero of $f$ at $0$, i.e.~the rate of convergence
of $f$ to zero, with the property that $\widehat{f}$ is in some
weighted $H^{p}$ space. See theorem 1.1 in \cite{MS02}. We omit
the details. This shows, for example, that $\widehat{f}(\xi)\neq O(e^{-\xi^{1/2+\epsilon}})$
for any $\epsilon>0$.

\end{rem*}
\begin{proof}
[Proof of theorem \ref{thm:existpoly}]Let $B$, $\tau'$ and $Y$
be the quadratic form, constant and lattice given for $R$ and $K+1$
by lemma \ref{lem:multidim}. Denoting\[
\tilde{a}(z):=\sum_{n:(n,z)\in Y}\tau'n^{-d/2}-P_{n}(z)\]
we see that $\tilde{a}(z)-a(z)$ is a constant (indeed, the precaution
of adding $\tau'n^{-d/2}$ is only necessary in the case $d=2$ in
order to make the sum converge). The fact that $R$ is bounded shows
that $\mathbb{P}(R_{n}=z)$ is zero for $n<c||z||$ and of course
$\sum_{n<c||z||}e^{-B(z)/n}=O(e^{-c||z||})$. Therefore we can use
(\ref{eq:reslem}) to get\[
\tilde{a}(z)=\sum_{n:(n,z)\in Y}\tau'n^{-d/2}\big(e^{-B(z)/n}Q_{2}(n,z)-1\big)+\sum_{n>c||z||}O(n^{-K-1})+O(e^{-c||z||})\]
Lemma \ref{lem:skhumnae} allows to calculate the left sum, while
the right one is obviously $O(||z||^{-K})$. This concludes the proof
of theorem \ref{thm:existpoly}.
\end{proof}
\label{constants}We see now that the quadratic form of theorem \ref{thm:existpoly}
is the same as that of lemma \ref{lem:multidim}, i.e.~$B(z)=\frac{1}{2}\langle M^{-1}z,z\rangle$
where $M$ is the correlation matrix of $R$. Thus $A=M^{-1/2}$ ---
of course, this is not uniquely defined, but the quadratic form $||Az||^{2}$
is. As for $\tau$, the factor $\vol Y\cap Z$ of lemma \ref{lem:multidim}
with the term $\frac{1}{q}$ of lemma \ref{lem:skhumnae} gives $\vol Z$
(this is easy to see) and we get \[
\tau=\frac{\vol Z}{(2\pi)^{d/2}|\det M|^{1/2}}\cdot\begin{cases}
2 & d=2\\
\Gamma(\frac{1}{2}d-1) & d\neq2\end{cases}\]
(the factor $2$ in dimension $2$ comes from the fact that we formulated
$\mu_{2}=\log||z||$ rather than $\log||z||^{2}$).

\section{The differential approach}

We shall demonstrate the differential approach by proving a few results
that would seem quite difficult using the direct approach only.

For simplicity of notation we shall always assume that the correlation
matrix of $R$ is constant. Thus we call a random walk \emph{spherical}
if it is bounded, balanced, and $\mathbb{E}\langle R,e_{i}\rangle\langle R,e_{j}\rangle=C\delta(i,j)$
(again, $e_{i}$ are unit vectors).

\begin{thm}
\label{thm:klhalf}Let $R$ be a $2$-dimensional spherical random
walk on $\mathbb{C}$. Then \begin{equation}
a(z)=\alpha\log|z|+\lambda+\real\left(\sum\alpha_{lk}\frac{z^{l}}{|z|^{k}}\right)+O(|z|^{-K})\label{eq:complexrep}\end{equation}
and $\alpha_{lk}=0$ whenever $l<\frac{1}{2}k$.
\end{thm}
The sum here is over a finite collection of $l\geq0$ and $k\geq0$
with the condition $-1\geq l-k\geq-K$.

\begin{proof}
The complex representation (\ref{eq:complexrep}) is nothing but a
restatement of theorem \ref{thm:existpoly}: writing $x=\frac{1}{2}(z+\bar{z})$
and $y=\frac{1}{2i}(z-\bar{z})$ we get a representation in the form
\[
\sum\alpha_{ijk}\frac{z^{i}\bar{z}^{j}}{|z|^{k}}\quad.\]
but of course if both $i\neq0$ and $j\neq0$ then we can cancel one
of them with the $|z|^{k}$ factor. The fact that the end result is
real means that $\alpha_{0jk}=\overline{\alpha_{j0k}}$ and we get
(\ref{eq:complexrep}). Thus what we need to prove is the fact that
$\alpha_{lk}=0$ whenever $l<\frac{k}{2}$. We change the representation
again, writing \begin{equation}
a(z)=\beta\log(z\bar{z})+\sum_{k+l>-K}\beta_{kl}z^{k}\bar{z}^{l}+O(|z|^{-K})\label{eq:azzbar}\end{equation}
and now we have to show that $\beta_{kl}=0$ whenever both $k$ and
$l$ are negative (of course, $\beta_{kl}=0$ when $k+l>0$). Now,
the first thing to notice is that we may differentiate (\ref{eq:azzbar})
formally. If $R$ moves to $v_{i}$ with probability $p_{i}$ then\begin{align*}
(\Delta f)(z) & =\sum p_{i}f(z+v_{i})-f(z)=\\
 & =\sum_{n=1}^{K-1}\frac{1}{n!}\sum_{i}p_{i}\frac{\partial^{n}}{\partial^{n}v_{i}}f(z)+O\left(\frac{1}{K!}\sum_{i}p_{i}\max_{\eta\in[0,1]}\left|\frac{\partial^{K}}{\partial^{K}v_{i}}f(z+\eta v_{i})\right|\right)\end{align*}
where we used the one-dimensional Taylor expansion for the functions
$f(z+tv_{i})$ for each $i$ (the constant in the $O(\cdot)$ above
is of course simply $1$). Denote the $n$th differential operator
in the sum by $D_{n}$ and denote the maximum inside the $O(\cdot)$
by $M_{K}$. Clearly, \begin{equation}
D_{n}z^{k}\bar{z}^{l}=O(|z|^{k+l-n})\qquad M_{n}z^{k}\bar{z}^{l}=O(|z|^{k+l-n})\quad.\label{eq:Dreddeg}\end{equation}
Now, $\Delta a=0$ with (\ref{eq:azzbar}) and (\ref{eq:Dreddeg})
show that\begin{equation}
\left(\sum_{n=1}^{K-1}D_{k}\right)\left(\beta\log z\bar{z}+\sum_{k+l>-K}\beta_{kl}z^{k}\bar{z}^{l}\right)=O(|z|^{-K})\quad.\label{eq:formldif}\end{equation}
This is exactly what we mean when we say that (\ref{eq:azzbar}) can
be differentiated formally.

Now to the actual calculation. Notice that $D_{1}=0$ (because $R$
is balanced) and that $D_{2}$ is, up to multiplication with a constant,
the continuous Laplacian, because $R$ is spherical. Therefore\begin{equation}
D_{2}(z^{k}\bar{z}^{l})=2\gamma klz^{k-1}\bar{z}^{l-1}\label{eq:delzzbar}\end{equation}
for some constant $\gamma\neq0$. Assume that we know $\beta_{kl}$
for all $k+l>-m$, and that $\beta_{kl}=0$ if both are negative.
We write\[
\left(\sum_{n=1}^{m}D_{n}\right)\left(\beta\log z\bar{z}+\sum_{k+l>-m}\beta_{kl}z^{k}\bar{z}^{l}\right)=\sum\delta_{kl}z^{k}\bar{z}^{l}\]
we get that $\delta_{kl}=0$ if both $k$ and $l$ are negative. Moreover,
(\ref{eq:formldif}) for $K>m+2$ gives equations for $\beta_{kl}$
with $k+l=-m$, namely $\beta_{kl}=\frac{1}{2\gamma kl}\delta_{k-1,l-1}$.
This shows inductively that $\beta_{kl}=0$ when they are both negative,
and the theorem is proven.
\end{proof}
The case $l=0$ in theorem \ref{thm:klhalf} has some potential uses
for precise estimates of hitting probabilities, so it seems worth
of special mention:

\begin{cor*}
For every spherical random walk on $\mathbb{R}^{2}$, the high-order
expansion of the harmonic potential $a$, considered as a function
of a continuous variable, satisfies the superpolynomial estimate \[
\int_{0}^{2\pi}a(Re^{it})\, dt=A\log R+B+O(R^{-K})\quad\forall K.\]

\end{cor*}
\begin{rem*}
\label{page:spherical}All the above discussion has analogues in dimension
$>2$, using the notions of spherical harmonics. Spherical harmonics
are functions on $S^{d-1}$ which are eigenvalues of the angular part
of the continuous Laplacian (which is the angular momentum squared
operator). It turns out that each is a restriction to $S^{d-1}$ of
a polynomial on $\mathbb{R}^{d}$ and we shall use this representation.
For a general introduction to the topic, see \cite[sect.\ 4.2]{SW71}.
\begin{itemize}
\item In two dimensions, a factor $P/|z|^{\alpha}$ (where $P$ is a homogeneous
polynomial) may be decomposed into a sum of the terms $z^{k}\bar{z}^{l}$,
$k+l=\deg P-\alpha$. The high dimension analogue is a unique decomposition
into a sum of terms $Q/||z||^{\beta}$ where $Q$ is a homogeneous
(continuously) harmonic polynomial, $\deg Q-\beta=\deg P-\alpha$. 
\item The analogue of the fact that only $z^{k}$ and $\bar{z}^{l}$ are
harmonic is the fact that a term $Q/||z||^{\beta}$ is harmonic if
and only if $\beta=2\deg Q+d-2$. This follows from the following
general formula:\begin{equation}
\Delta_{C}\left(\frac{Q}{||z||^{\beta}}\right)=\frac{||z||^{2}\Delta_{C}Q+(\beta^{2}-\beta(2\deg Q+d-2))Q}{||z||^{\beta+2}}\label{eq:LapQzalpha}\end{equation}
and from the harmonicity of $Q$ (remember that $\Delta_{C}$ is the
continuous Laplacian). (\ref{eq:LapQzalpha}) also gives the analogue
of (\ref{eq:delzzbar}).
\item The analogue of {}``$k$ and $l$ are both negative'' is {}``$\beta<2\deg Q+d-2$'',
or, equivalently, {}``$\Delta_{C}^{n}(Q/||z||^{\beta})$ is never
zero, for any $n$''.
\item The analogue of theorem \ref{thm:klhalf} is the following: if $R$
is a spherical random walk and\[
a(z)=\frac{\tau}{||z||^{d-2}}+\lambda+\frac{Q(z)}{||z||^{\beta}}+O(||z||^{-K})\]
then\[
\Delta_{C}^{K}\left(\frac{Q(z)}{||z||^{\beta}}\right)=0\]
We omit the proof, which is a similar induction.
\item The corollary to theorem \ref{thm:klhalf} is also true in higher
dimensions. This follows from the fact that the terms $Q/||z||^{\beta}$
for different $\beta$'s are orthogonal in the sphere measure: \[
\int_{S^{d-1}}Q_{1}\cdot Q_{2}=0\quad\forall Q_{1},Q_{2}\textrm{ harmonic, }\deg Q_{1}\neq\deg Q_{2}\]
and in particular $Q$ is orthogonal to $1$ if $\deg Q>0$.
\end{itemize}
\end{rem*}
\begin{thm}
\label{thm:eqssuff}The coefficients of the high-order expansion of
the harmonic potential $a$ of the regular random walk on $\mathbb{Z}^{2}$
can be calculated from the differential equations and (\ref{eq:anexplic}).
\end{thm}
\begin{proof}
There is almost nothing to prove here. We saw during the proof of
the previous theorem that (\ref{eq:formldif}) allows to calculate
all coefficients except those where $k=0$ or $l=0$ (i.e.~the coefficients
of the harmonic summands). There is only one of these (and its conjugate)
in every level $k+l=-n$. Therefore the coefficient of $z^{-n}$ can
be determined from the coefficients of $z^{-n-1}\bar{z},z^{-n-2}\bar{z}^{2},\dotsc$
and from (\ref{eq:anexplic}).
\end{proof}
The proof of theorem \ref{thm:eqssuff} is also the second algorithm
for getting (\ref{eq:a_high_order}).

\begin{thm}
\label{thm:deg} \hfill



\begin{enumerate}
\item The term of order $-k$ in any walk is $P/||z||^{\alpha}$ where $\alpha\leq4k+2-d$
and $P$ is a polynomial of degree $\alpha-k\leq3k+2-d$.
\item If $D_{3}\equiv0$, in particular if $R$ is reversible%
\footnote{A reversible (or symmetric) walk is a random walk satisfying $R=-R$.
In graph language, a reversible random walk is a walk on a simple
graph, while a non-reversible walk is a walk on a directed graph.
Many of the claims of this paper need slight alterations for a non-reversible
walk. For example, we have $\tilde{\Delta}a=\delta_{\{0\}}$ where
$\tilde{\Delta}$ is the Laplacian of the time-reversed walk.%
}, then the term of order $-k$ is $P/||z||^{\alpha}$ where $\alpha\leq3k+2-d$
and $P$ is a polynomial of degree $\alpha-k\leq2k+2-d$.
\end{enumerate}
\end{thm}
\begin{proof}
Take as an example clause 2. We may assume that the walk is spherical.
In two dimensions, the theorem is equivalent to showing that only
terms of the form $z^{k}\bar{z}^{l}$, $0\leq k\leq\frac{1}{3}|l|$
(or $0\leq l\leq\frac{1}{3}|k|$) appear. On the one hand, $k$ can
only increase (by $1$) in the step of taking $D_{2}^{-1}$ (see the
proof of theorem \ref{thm:klhalf}), and this operation also increases
$l$ by $1$. On the other hand, the step of taking $D_{4}+D_{5}+\dotsb$
decreases $l$ by at least $4$, and the theorem follows. In higher
dimensions the same argument works using the decomposition into spherical
harmonics discussed before theorem \ref{thm:eqssuff}. 
\end{proof}
We note that clause 1 of theorem \ref{thm:deg} is already mentioned
in \cite{FU96,U98}.

Theorem \ref{thm:deg} can be generalized to cases where higher symmetries
exist. For example, for the $6$-walk on the triangular lattice, $D_{3}\equiv D_{5}\equiv0$
and $D_{4}=\frac{9}{4}\Delta_{C}^{2}$ which allow to get that the
term of order $|z|^{-k}$ is $P/|z|^{\nu}$ where $\nu\leq\frac{5}{2}k$
and $\deg P\leq\frac{3}{2}k$.

\begin{thm}
The harmonic potential of a reversible random walk contains only terms
of even order in even dimension and only terms of odd order in odd
dimension.
\end{thm}
\begin{proof}
Since $D_{2k+1}\equiv0$ for all $k$, we only need to handle the
terms we cannot calculate, namely the (continuously) harmonic terms.
However, a harmonic term of these orders must be $P/||z||^{\alpha}$
with $\deg P$ odd, and since $a(z)=a(-z)$ such terms cannot appear
in the expansion of $a$.
\end{proof}

\section{\label{sec:Explicit-constants}Explicit constants}

This note is a spin-off from a different work \cite{K} where it was
necessary to have an \emph{explicit} value for the constants in the
$O(\cdot)$ in (\ref{eq:stoher}). How should one go about to calculate
something like that? Generally, one would try to get some explicit
constant in the $O$ of an approximation of one order higher, and
then use the McCrea-Whipple algorithm to check the first few values.

Which approach would be best to get an explicit constant in the $O(\cdot)$?
The differential approach doesn't seem to give any explicit constants
whatsoever. Therefore it was necessary to use the direct approach.
A program was written to handle all the algebra and the calculations
of the explicit constants in the $O(\cdot)$.%
\footnote{The program is several thousands of lines long. It is available from
the first author for examination, further development etc.%
} After careful selection of parameters, the program got (\ref{eq:a_high_order})
with an error term which is smaller than $\leq10^{25}|z|^{-12}$ for
any $|z|>1200$. The program doing the algebra ran for half an hour
on a PC. This may or may not be enough --- it depends on the difference
between the actual maximum and the asymptotic value. However, this
difference is $>0.01$ (see below) and $10^{25}|z|^{-12}<10^{-5}|z|^{-2}$
for $|z|>1200$ so this quantity is indeed negligible.

Interestingly, here the third approach (the one discussed briefly
on page \pageref{Thirdapproach}) shows a slight advantage over the
other two. The third approach allows to get an explicit constant in
the $O(\cdot)$ in (\ref{eq:thirdstoher}) with much less effort ---
a few trivial programs to find optimal parameters were all that was
necessary. Contrariwise to the simplicity of the programs, the estimates
are not as good and it was necessary to run the McCrea-Whipple algorithm
for $|z|<6000$, which takes a few days.

The maximum of $|z|^{2}\left(a(z)-\frac{2}{\pi}\log|z|-\lambda\right)$
was found at $z=3$ and so the constant implicit in the $O(\cdot)$
in (\ref{eq:stoher}) is \[
\frac{432+9\log72+18\gamma}{\pi}-153\approx0.06882\]
Notice that this is only slightly larger than the asymptotics $\frac{1}{6\pi}\approx0.05305$.
Naturally, explicit constants can be calculated for higher order estimates
as well.

\subsection{Other simulations}

We add below two additional $a$'s. The first is for a {}``drunk
king's walk'' --- a random walk on $\ZZ^{2}$ going at each step
to each of the 8 neighboring lattice points with equal probability.
Here are the first few coefficients:\begin{align*}
a(z) & =\frac{4}{3\pi}\log|z|+\lambda+\real\left(\frac{z^{4}}{9\pi|z|^{6}}+\frac{11z^{4}}{90\pi|z|^{8}}-\frac{5z^{8}}{36\pi|z|^{12}}-\frac{167z^{8}}{252\pi|z|^{14}}+\right.\\
 & \left.+\frac{35z^{12}}{54\pi|z|^{18}}-\frac{1673z^{8}}{2160\pi|z|^{16}}+\frac{15z^{12}}{2\pi|z|^{20}}-\frac{1925z^{16}}{288\pi|z|^{24}}+\dotsb\right)\end{align*}
 so in this case it is no longer true that all the coefficients have
a constant sign. Another interesting case is the non-reversible walk
on the directed triangular lattice going at each step with probability
$\frac{1}{3}$ to $z+1$, $z+\omega$ or $z+\omega^{2}$ where $\omega=-\frac{1}{2}+i\frac{\sqrt{3}}{2}=\sqrt[3]{1}$.
Here are the first few coefficients:\begin{align*}
a(z) & =\frac{\sqrt{3}}{\pi}\log|z|+\lambda+\frac{\sqrt{3}}{\pi}\real\left(\frac{z^{3}}{6|z|^{4}}+\frac{z^{6}}{12|z|^{8}}+\frac{z^{3}}{18|z|^{6}}+\frac{5z^{9}}{54|z|^{12}}+\frac{17z^{6}}{135|z|^{10}}+\right.\\
 & \left.\frac{35z^{12}}{216|z|^{16}}+\frac{19z^{9}}{54|z|^{14}}+\frac{7z^{15}}{18|z|^{20}}+\frac{85z^{6}}{1134|z|^{12}}+\frac{98z^{12}}{81|z|^{18}}+\frac{385z^{18}}{324|z|^{24}}+\dotsb\right)\end{align*}
We see that the phenomenon of constant sign is not specific to the
regular random walk (we checked this up to coefficients of order $O(|z|^{-28})$).
Do note that in this case all coefficients are positive while in the
$\mathbb{Z}^{2}$ case all coefficients are negative.

\end{document}